\documentclass[12pt]{amsart}

\usepackage{mathrsfs}  
\usepackage{amssymb}
\usepackage{verbatim}

\usepackage[all]{xy}
\usepackage{pb-diagram,pb-xy}
\title[Homotopy type of the Deligne--Mumford compactification]{On the homotopy type of \\the Deligne--Mumford compactification}

\author{Johannes Ebert, Jeffrey Giansiracusa} 

\theoremstyle{theorem}
\newtheorem{theorem}{Theorem}[section] 
\newtheorem{lemma}[theorem]{Lemma} 
\newtheorem{proposition}[theorem]{Proposition}
\newtheorem{corollary}[theorem]{Corollary}
\theoremstyle{remark}

\newtheorem{remark}[theorem]{Remark}

\newcommand{\Obj}{\mathrm{Obj}}
\newcommand{\Mor}{\mathrm{Mor}}

\newcommand{\Ho}{\mathrm{Ho}}
\newcommand{\bC}{\mathbb{C}}
\newcommand{\C}{\mathbb{C}}
\newcommand{\Q}{\mathbb{Q}}
\newcommand{\bH}{\mathbb{H}}
\newcommand{\bR}{\mathbb{R}}

\newcommand{\Z}{\mathbb{Z}}
\newcommand{\bZ}{\mathbb{Z}}
\newcommand{\M}{\mathcal{M}_{g,n}}
\newcommand{\bM}{\overline{\mathcal{M}}_{g,n}}
\newcommand{\bMan}{\overline{\mathcal{M}}\hbox{}_{g,n}^{an}}
\newcommand{\MCG}{\mathcal{MCG}_{g,n}}
\newcommand{\fB}{\mathfrak{D}}
\newcommand{\CL}{\mathcal{CL}}

\newcommand{\co}{:}

\begin{document}
\begin{abstract}
  An old theorem of Charney and Lee says that the classifying space of
  the category of stable nodal topological surfaces and isotopy
  classes of degenerations has the same rational homology as the
  Deligne--Mumford compactification.  We give an integral refinement:
  the classifying space of the Charney-Lee category actually has the
  same homotopy type as the moduli stack of stable curves, and the
  \'etale homotopy type of the moduli stack is equivalent to the
  profinite completion of the classifying space of the Charney-Lee category.
\end{abstract}
\maketitle

\section{Introduction}

The purpose of this note is to give an integral refinement of a
relatively old theorem of Charney and Lee \cite{Charney-Lee} giving a
model for the rational homology of the Deligne--Mumford
compactification of the moduli space of curves in terms of a category
made of mapping class groups.

Let $\M$ denote the moduli stack of proper smooth algebraic curves of
genus $g$ with $n$ ordered marked points, and let $\bM$ denote the
moduli stack of stable curves (the Deligne--Mumford compactification of
$\M$). They are both smooth Deligne--Mumford algebraic stacks defined
over $\mathrm{spec} \: \Z$.  These algebraic stacks have associated
complex analytic stacks (orbifolds), $\M^{an}$ and $\bMan$.  It is
well known that the coarse moduli space of $\M^{an}$ has the same
rational homology as the classifying space of the mapping class group
$\MCG$ of a surface of genus $g$ with $n$ marked points.

Charney and Lee defined a category $\CL_{g,n}$ in which:
\begin{itemize}
\item objects are stable nodal surfaces of genus $g$ with $n$ ordered
  distinct marked points in the smooth part,
\item morphisms are isotopy classes of orientation-preserving
  diffeomorphisms and degenerations (a degeneration is a map which
  collapses some circles to nodes and is a diffeomorphism on the
  complement of these circles) that respect the marked points.
\end{itemize}
The mapping class group $\MCG$ sits inside $\CL_{g,n}$ as the
automorphism group of a smooth surface; automorphism groups of other
objects are mapping class groups of singular surfaces appearing in the
boundary of the Deligne--Mumford compactification.  Note that the
moduli stack $\bM$ and the category $\CL_{g,n}$ both have
stratifications by `dual graphs'.

Charney and Lee proved \cite[Theorem 6.1.1]{Charney-Lee} that (for
$n=0$) the classifying space of $\CL_g$ has the same rational homology
as the coarse moduli space of $\mathcal{M}_{g}^{an}$.  The moduli
stack and the coarse moduli space have the same rational homology, but
integrally they differ!  The mod $p$ homology of the open moduli
\emph{stack} has been computed in the Harer-Ivanov stable range by
Galatius \cite{Galatius} (using the theorem of Madsen and Weiss
\cite{Madsen-Weiss}); it contains much more than just reductions of
non-torsion classes.  The mod $p$ homology of the Deligne--Mumford
compactified stack has been studied by Galatius--Eliashberg
\cite{Galatius-Eliashberg} and the authors
\cite{Ebert-Giansiracusa-PT}, but it remains largely unknown.

An analytic stack (or more generally a topological stack)
$\mathfrak{X}$ has a homotopy type which can be defined by choosing a
covering $X \to \mathfrak{X}$ by a space $X$ and then taking the
geometric realization of the simplicial space which in degree $n$ is
the $(n+1)$-fold fiber-product $X \times_{\mathfrak{X}} \cdots
\times_{\mathfrak{X}} X$.  I.e. take the classifying space of a
topological groupoid presenting the orbifold; see Moerdijk
\cite{Moerdijk}, Noohi \cite{Noohi}, Ebert--Giansiracusa
\cite{Ebert-Giansiracusa-PT} and Ebert \cite{Ebert} for more details.  
The integral singular homology and fundamental group of the
analytic stack agree with those of the homotopy type.  As an example
of homotopy types, it is well know that the homotopy type of the stack
$\M^{an}$ is $B\MCG$

We prove the following integral refinement of Charney and Lee's theorem.

\begin{theorem}\label{main-theorem}
  The classifying space of $\CL_{g,n}$ is homotopy equivalent to the
  homotopy type of the stack $\bMan$, so in particular, $H^*(\bMan;\Z)
  \cong H^*(B\CL_{g,n};\Z)$.  Furthermore, this homotopy equivalence
  is compatible with the stratifications of $\CL_{g,n}$ and $\bMan$
\end{theorem}

By Artin--Mazur \cite{Artin-Mazur}, Oda \cite{Oda} and
Frediani--Neumann \cite{Neumann} a Deligne--Mumford algebraic stack
has an \'etale homotopy type (living in the category of pro-objects in
the homotopy category of simplicial sets).  By the Comparison Theorem
of \'etale homotopy theory \cite[Theorem 8.4]{Friedlander}, the
\'etale homotopy type of a stack over $\overline{\Q}$ is weakly
equivalent to the Artin--Mazur profinite completion of
the homotopy type of the associated analytic stack.  Let $\M \otimes
\overline{\Q}$ denote the extension of scalars of $\M$ to
$\overline{\Q}$ (i.e. the restriction of the moduli functor to schemes
over $\overline{\Q}$).  As explained in \cite{Oda}, the \'etale
homotopy type of $\M \otimes \overline{\Q}$ is the Artin--Mazur
profinite completion $(B\MCG)^{\wedge}$.  Similarly, Frediani--Neumann
\cite{Neumann} describes the \'etale homotopy type of the moduli stack
of curves with an action of a finite group $G \subset \MCG$.  In this
vein, the Comparison Theorem plus Theorem \ref{main-theorem} yields:
\begin{corollary}
  The \'etale homotopy type of $\bM\otimes \overline{\Q}$ is weakly
  equivalent to the Artin--Mazur profinite completion
  $(B\CL_{g,n})^{\wedge}$, and this equivalence respects the respective
  stratifications.
\end{corollary}
(Recall that a weak equivalence of pro-objects is a morphism of
pro-objects which induces an isomorphism on their homotopy pro-groups.)

The original Charney-Lee proof could probably easily be adapted to
handle surfaces with marked points and to show that the rational
homology equivalence is compatible with the stratifications.  However,
our proof is significantly more direct than theirs, while also giving
the integral refinement.  Our proof is based on existence of a
particularly nice atlas, first constructed by Bers \cite{Bers1, Bers2,
  Bers3, Bers4}, which is well-adapted to the combinatorial structure
of the stratification of $\bM$.  Roughly speaking, the Bers atlas
generalizes the Teichm\"uller space in the same way that the
Charney-Lee category generalizes the mapping class group.

\subsection*{Acknowledgements}

The first author was supported by a postdoctoral grant of the German
Academic Exchange Service (DAAD). He enjoyed the hospitality of the
Mathematical Institute of the University of Oxford.  The second author
thanks C.-F. B\"odigheimer for an invitation to visit Universit\"at
Bonn in April 2006, where this project was begun.  The project was
completed during the second author's stay at the IHES during
2006/2007, and he thanks both institutions for their hospitality.

\section{The Charney-Lee category}

Before proceeding with the principal content of this note, we collect
here some remarks about the Charney-Lee category.  We will not need
either of these remarks, so we only sketch them briefly, but the
reader might nevertheless find these comments illuminating.

Firstly, there is a topological version $\CL_{g,n}^{top}$; is has the
same objects as $\CL_{g,n}$, while the space of morphisms
$\CL_{g,n}^{top}(S,T)$ is the space of degeneration maps $S\to T$
(i.e. maps which collapse circles to nodes and which are
orientation-preserving diffeomorphisms outside these collapsed
circles).  The topology of the morphism spaces is the Whitney
$C^{\infty}$-topology.  We can clearly identify
$\pi_0(\CL_{g,n}^{top})$ with $ \CL_{g,n}$.  Moreover, the obvious
functor $\CL_{g,n}^{top} \to \CL_{g,n}$ is a homotopy equivalence of
categories; in other words, the components of the morphism spaces in
$\CL_{g,n}^{top}$ are all contractible.  This is a generalization of
the well-known theorem \cite{EaEe, EaSch, Gram} that the
diffeomorphism groups of oriented smooth surfaces with boundary (with
negative Euler number) have contractible components.  The space of
degenerations $\CL^{top}_{g,n}(S,T)$ fibers over the space of
unparametrized 1-dimensional submanifolds in $S$ by
taking the union of all curves which are collapsed.  This map is a
Serre fibration and the fibers are homeomorphic to diffeomorphism
groups of surfaces with negative Euler characteristic, hence the
components of the fibers are contractible.  It follows from
\cite{Epst,Gram} that the components of the base are also
contractible.

A second remark is that $\CL_{g,n}$ can be described as an orbit
category.  The \emph{orbit category} of $\MCG$ is the category whose
objects are orbits $\MCG / H$ and whose morphisms are the
$\MCG$-equivariant maps.  The category $\CL_{g,n}$ is equivalent to
the full subcategory of the orbit category containing precisely
those orbits for which the isotropy subgroup $H$ is a free abelian
group generated by a collection of disjoint Dehn twists.  To see this,
fix a smooth surface $S$ of genus $g$ with $n$ marked points, and for
each object $T \in \CL_{g,n}$ choose a degeneration $p(T)\co S \to T$.
The Dehn twists on $S$ around the inverse images of the nodes of $T$
determine a free abelian subgroup of $\MCG$ and hence an orbit $O(T)$.
Given a degeneration $\alpha\co T \to T'$, there exists
$\widetilde{\alpha} \in \MCG$ such that $\alpha \circ p(T) = p(T')
\circ \widetilde{\alpha}$ --- this $\widetilde{\alpha}$ is only unique
up to certain Dehn twists, but it induces a well-defined morphism
$O(T) \to O(T')$ in the orbit category.

\section{The Bers atlas for $\bMan$}\label{bers-groupoid}

Bers \cite{Bers1, Bers2, Bers3, Bers4} has constructed an atlas $\fB$,
which we shall call the \emph{Bers atlas}, for the differentiable
stack $\bMan$.  (To avoid notational clutter we leave $g$ and $n$
implicit).  This atlas is an extension of the atlas for the
uncompactified moduli stack $\M^{an}$ given by Teichm\"uller space.

The Bers atlas is defined as follows.  Let $S$ be a fixed stable nodal
topological surface of genus $g$ with $n$ marked points.  An
\emph{$S$-marked Riemann surface} is a stable nodal Riemann surface
$F$ with $n$ marked points lying in the smooth part, together with a
degeneration $F \to S$ which respects the marked points.  Two
$S$-marked Riemann surfaces $f\co F\to S$ and $f'\co F' \to S$ are defined
to be equivalent if there exists a biholomorphic map $g\co F
\stackrel{\cong}{\to} F'$ (respecting the marked points) such that the
diagram
\[
\begin{diagram}
\node{F} \arrow{se,b}{f} \arrow[2]{e,t}{g} \node[2]{F'} \arrow{sw,r}{f'} \\
\node[2]{S}
\end{diagram}
\]
commutes up to a homotopy that is constant on the marked points.

Let $\fB(S)$ denote the set of all equivalence classes of $S$-marked
Riemann surfaces.  In \cite{Bers2, Bers4} Bers defined a topology on
$\fB(S)$ making it into a \emph{contractible} manifold, and such that
when $S$ is smooth then $\fB(S)$ is the usual Teichm\"uller space of
$S$.

In fact, the Fenchel-Nielsen coordinates give a homeomorphism between
$\fB(S)$ and an open ball as follows.  Let $N$ denote the set of nodes
of $S$ and choose a complete cutsystem $C$ on $S$ (i.e. a collection
of disjoint simple closed curves in the smooth part of $S$ such that
the complement of $C \sqcup N$ is a disjoint union of pairs of pants.
Given a point $[f\co F\to S] \in \fB(S)$, there is a unique compatible
hyperbolic metric on $F$.  The free homotopy class of each curve of
$f^{-1}(C)$ has a minimal geodesic length and a twist; these numbers determine
a point in $(\bR_+ \times \bR)^C \cong \bH^C$.  For a node $n\in N$,
if $f^{-1}(n)$ is a simple closed curve then this free homotopy class
has a length and a twist in $\bR_+\times \bR/\bZ \cong \bC^\times$,
and the coordinates converge to the origin as the inverse image of $n$
in $F$ collapses to a node.  Hence the Fenchel-Nielsen coordinates
give a map
\[
\fB(S) \to \bH^C \times \bC^N,
\]
which one can show is a homeomorphism.  A particularly nice exposition
for smooth surfaces can be found in Hubbard's book \cite[p. 320 ff]{Hubb}. 

Bers also endowed $\fB(S)$ with the structure of a complex manifold
which embeds as a bounded domain in $\C^{3g - 3 + n}$, generalizing
the Maskit coordinates, but we shall not need this fact.

The Bers atlas is given by
\[
\coprod_S\fB(S) \to \bMan,
\]
where the disjoint union runs over each diffeomorphism class of stable
nodal surfaces $S$ having genus $g$ and $n$ marked points; the map to
the moduli space is given informally by forgetting the markings,
sending a marked Riemann surface $[F \to S]$ to $F$.  More precisely,
there is a tautological family over $\fB(S)$ whose fiber over $[F\to
S]$ is $F$, and the map to $\bMan$ is given by classifying this
tautological family.

\begin{theorem}
  The morphism $\coprod_S\fB(S) \to \bMan$ defines a proper \'etale atlas for $\bMan$
  as a differentiable (or even complex analytic) stack.
\end{theorem}
\begin{proof}
This is essentially contained in the work of Bers; it follows from
Theorems 6 and 7 announced in \cite{Bers2}.
\end{proof}

Put differently, the representable submersion $\coprod_S \fB(S) \to
\bMan$ determines a Lie groupoid which is proper and \'etale (i.e. an
orbifold groupoid, though not always effective),
\[
\fB := \left[ \left( \coprod_S\fB(S) \right)\times_{\bMan}
\left(\coprod_S\fB(S) \right) \rightrightarrows
\left(\coprod_S \fB(S) \right) \right];
\]
we call this the \emph{Bers Groupoid}.  An object of this groupoid is
the equivalence class of an $S$-marked Riemann surface $F$ for some
$S$; a morphism $[F\to S] \to [F'\to T]$ is a biholomorphic map
$g\co F \stackrel{\cong}{\to} F'$ respecting the marked points but
completely ignoring the maps to $S$ and $T$.  We call this the Bers
groupoid and denote it $\fB$.  Since it is a presentation of the stack
$\bMan$, its classifying space is a model for the homotopy type of
$\bMan$.  In particular, $B\fB$ has the same integral (co)homology as
the stack $\bMan$.

We now recall some facts about the Bers atlas from \cite{Bers2,Bers4}.
A degeneration $\alpha\co S \to T$ induces a map $\alpha_*\co \fB(S) \to
\fB(T)$ by change-of-marking, i.e.
\[
[F\stackrel{f}{\to} S] \mapsto [F \stackrel{\alpha\circ f}{\to} T].
\]
The induced map $\alpha_*$ is a local homeomorpism.  Its image is
precisely the subspace consisting of those points $[F\to T]$ for which
the marking can be lifted along $\alpha$ to $S$; with appropriate
Fenchel-Nielsen coordinates one easily sees that this is the
complement of a collection of complex coordinate hyperplanes.  The map
$\alpha_*$ only depends on the isotopy class of $\alpha$ because of
the equivalence relation on degenerations $F\to S$ used in defining
the space $\fB(S)$. In particular, there is a properly discontinuous
action of the mapping class group $\mathcal{MCG}(S)$ of $S$ on
$\fB(S)$ and the quotient stack $[\fB(S)/\mathcal{MCG}(S)]$ is
isomorphic to the image of $\fB(S)$ in $\bMan$.

\section{The Bers groupoid and the Charney-Lee category}

We shall now describe a subcategory of the Bers groupoid which is more
visibly related to the Charney-Lee category.  We then give a
completely explicit description of the Bers groupoid in terms of this
subcategory.

The spaces $\fB(-)$ together with the change-of-marking maps described
above determine a functor $\widehat{\fB}\co \CL_{g,n} \to
\mathrm{Spaces}$, and we may form the transport category (or
Grothendieck construction) $\CL_{g,n} \int \widehat{\fB}$.
Concretely, an object of the transport category is a point $[f\co F\to
S]$ in $\fB(S)$ for some $S$.  A morphism from $[f\co F\to S]$ to $[f'\co
F'\to T]$ is represented \emph{a priori} by a biholomorphic map $g\co F
\to F'$ together with the isotopy class of a degeneration $\alpha\co S
\to T$ such that the diagram
\[
\begin{diagram}
\node{F} \arrow{e,t}{g} \arrow{s,l}{f} \node{F'} \arrow{s,r}{f'} \\
\node{S} \arrow{e,t}{\alpha} \node{T}
\end{diagram}
\]
commutes up to homotopy. However, since the Charney-Lee category
possesses the right cancelation property,
\[
[\alpha] \circ [\gamma] = [\beta] \circ [\gamma] \mbox{\: implies \:}
[\alpha] = [\beta],
\]
the isomorphism $g$ uniquely determines the degeneration isotopy class
$[\alpha]$.  Note that not every isomorphism covers (up to homotopy) a
degeneration.  Thus a morphism $[F\to S] \to [F' \to T]$ can be
specified simply by a biholomorphic map $g\co F \to F'$ for which
there exists degeneration isotopy class that it covers.

By comparing the definitions the following is now apparent.
\begin{proposition}\label{inclusion-prop}
  The topological category $\CL_{g,n} \int \widehat{\fB}$ is
  isomorphic to a subcategory of the Bers groupoid $\fB$; namely, it is the
  subcategory with all objects of $\fB$ and only those biholomorphic
  maps which cover (up to isotopy) a degeneration of the markings.
\end{proposition}

\begin{lemma}\label{inclusion-htpy-equiv-lemma}
The inclusion $\CL_{g,n} \int \widehat{\fB} \hookrightarrow \fB$ induces a
homotopy equivalence of classifying spaces.
\end{lemma}

We will give the proof of Lemma \ref{inclusion-htpy-equiv-lemma} in
section \ref{proof-of-main-lemma} after some preparation in section
\ref{technicalities}.  Assuming this lemma for the moment, the proof
of Theorem \ref{main-theorem} is straightforward.

\begin{proof}[Proof of Theorem \ref{main-theorem}]
The Bers atlas is an atlas for $\bMan$, and so by definition of the
homotopy type of a stack, the classifying space $B\fB$ of the Bers
groupoid $\fB$ is the homotopy type of $\bMan$. Because
$\widehat{\fB}$ takes any $S \in \CL_{g,n}$ to a contractible space,
the forgetful functor $\CL_{g,n} \int \widehat{\fB} \to \CL_{g,n}$
induces a homotopy equivalence on classifying spaces. Therefore, by
Lemma \ref{inclusion-htpy-equiv-lemma}:
\[
 \Ho(\bMan) \simeq  B\fB  \stackrel{\simeq}{\leftarrow} 
 B\left({\CL_{g,n}\int \widehat{\fB}}\right) \stackrel{\simeq}{\to} B \CL_{g,n}.
\]
We postpone the discussion of compatibility with the stratifications
until section \ref{Stratifications}.
\end{proof}

\begin{remark}
  It is possible to show that, as abstract categories, when one
  formally adjoins inverses to all arrows of $\CL_{g,n}\int
  \widehat{\fB}$ then one obtains precisely the Bers groupoid $\fB$.  In
  particular, an arrow $[f\co F \to S] \to [f'\co F' \to T]$ can be
  represented by $\alpha^{-1} \circ \beta$ for a pair of degenerations
  $S\stackrel{\alpha}{\leftarrow} F' \stackrel{\beta}{\to} T$, and
  this representation is unique up to precomposition with an element
  of the mapping class group of $F'$.
\end{remark}

\section{A lifting property of the Bers atlas}\label{technicalities}

Let $X$ be a space and $\sigma\co X \to \bMan$ be a map. We say that a
lift $\widetilde{\sigma}\co X\to \fB(S)$ of $\sigma$ is \emph{maximal}
if $\widetilde{\sigma}$ does not admit a lift to $\fB(S')$ for any
$S'$ with a strict degeneration $S' \to S$.  Clearly, if $\sigma$
admits a lift to some $\fB(T)$ then it lifts further to a maximal
lift.

The goal of the present section is to prove the following result.

\begin{lemma}\label{initial-lifting-lemma}
Suppose $X$ is simply connected and $X \to \bMan$ admits maximal
lifts $\sigma_1\co X \to \fB(S)$ and $\sigma_2\co X \to \fB(T)$.  Then
there exists a diffeomorphism (unique up to isotopy) $\alpha\co S \cong
T$ with $\alpha_* \sigma_1 = \sigma_2$.
\end{lemma}

An equivalent formulation of the above lemma is that for any pair of
stable surfaces $S$ and $T$, there exists a stable surface $R$
degenerating onto $S$ and $T$ such that the map from $\fB(R)$ to any
component of the universal cover of $\fB(S) \times_{\bMan} \fB(T)$ is
a homeomorphism.  However, we do not know a more direct proof of
this fact.

The main tool for the proof of Lemma \ref{initial-lifting-lemma} is a sheaf
of sets $\mathcal{Z}$ on the the differentiable stack $\bMan$.  This
sheaf encodes the continuity property of markings on the fibers in a
marked family of stable Riemann surfaces.  The idea of the sheaf is as
follows.  Given a family $E \to X$ of stable Riemann surfaces, an
element of $\mathcal{Z}(X)$ should be thought of as the isotopy class
of a continuous subfamily $C \subset E$ that restricts in each fiber
to either a node or a simple closed curve that does not not meet the
nodes and marked points and does not retract to a node.  If $X\to
\fB(S)$ is a maximal lift then one can reconstruct the homeomorphism
type of $S$ from the sections of $\mathcal{Z}$ over $X$ that restrict
to a node in some fiber: each node in each fiber determines a node of
$S$ and maximality of the lift ensures that $S$ has no superfluous
nodes.  This will show that maximal lifts are essentially unique.

We construct the sheaf precisely by defining its stalks over the Bers
groupoid, topologizing its \'etale space, and then showing it descends
to a sheaf on the stack $\bMan$.  The stalk $\mathcal{Z}_{[S]}$ at a
point $[S] \in \bMan$ is defined to be the union of the set of nodes
of $S$ with the set of isotopy classes of unoriented simple closed
curves in $S\smallsetminus \{\mbox{nodes and marked points}\}$ which
bound neither a disc nor a once-punctured disc.  A degeneration
$\alpha\co T\to S$ induces an injective map $\alpha^*\co \mathcal{Z}_{[S]}
\hookrightarrow \mathcal{Z}_{[T]}$ by taking preimages of curves and
nodes.  Thus over $\fB(S)$ the markings canonically identify
$\mathcal{Z}_{[S]}$ with a subset of each stalk.

As a set, the \'etale space $\mathcal{Z}_{et}$ of $\mathcal{Z}$ is
the disjoint union of the stalks; we topologize it as follows. Given
a point $[\alpha\co F \to S] \in \fB(S)$, there exists a neighborhood
$U$ of this point which lifts along the local homeomorphism
$\alpha_*\co \fB(F) \to \fB(S)$, and so the markings identify
$\mathcal{Z}_{[S]}$ with a subset of each stalk over $U$. The
topology is determined by the condition that a section over $U$ is
continuous at $[F\to S]$ if and only if it is locally constant with
respect to these identifications.

Next we claim that $\mathcal{Z}$ is a sheaf on the stack $\bMan$. To
justify this, we need to argue that $\mathcal{Z}$ satisfies the
appropriate descent conditions. More precisely, let $d_0, d_1\co \Mor
(\fB) \to \Obj(\fB)$ be source and target maps and let $d_0, d_1,
d_2\co \Mor(\fB) \times_{\Obj(\fB)} \Mor(\fB) \to \Mor (\fB)$ be the
three simplicial structure maps in the nerve of $\fB$ (they are left
projection, composition, and right projection respectively).  A
descent datum for $\mathcal{Z}$ is an isomorphism $f\co d_{0}^{*}
\mathcal{Z} \to d_{1}^{*} \mathcal{Z}$ which makes the hexagon of
sheafs and isomorphisms on $\Mor(\fB) \times_{\Obj(\fB)} \Mor(\fB)$
commutative (the equalities are induced from simplicial identities):
\[
\xymatrix{
&  d_{0}^{*} d_{0}^{*} \mathcal{Z} \ar[r]^{d_{0}^{*} f}
& d_{0}^{*} d_{1}^{*} \mathcal{Z}   \ar@{=}[dr] &    \\
d_{1}^{*} d_{0}^{*} \mathcal{Z} \ar@{=}[ur] \ar[dr]^{d_{1}^{*} f} &
& & d_{2}^{*}
d_{0}^{*} \mathcal{Z} \ar[dl]^{d_{2}^{*} f} \\
&  d_{1}^{*} d_{1}^{*} \mathcal{Z} \ar@{=}[r]  & d_{2}^{*} d_{1}^{*}
\mathcal{Z} &
}
\]
There is an obvious bijection of \'etale spaces
$d_{0}^{*}\mathcal{Z}_{et} \cong d_{1}^{*} \mathcal{Z}_{et}$ and the
topology is designed so that this is a homeomorphism. The
commutativity of the above diagram is also clear. Thus $\mathcal{Z}$
is a sheaf of sets on $\bMan$.

The following property follows immediately from the definition of the
topology on the \'etale space of $\mathcal{Z}$.
\begin{lemma}\label{central-stalk-lemma}
  Let $U\subset \fB(S)$ be a neighborhood of the origin $[S \to S]$.
  There is a canonical bijection $\mathcal{Z}(U) \cong
  \mathcal{Z}_{[S]}$ induced in one direction by restriction to the
  stalk over $[S\to S]$ and in the other direction by using the
  markings to identify $\mathcal{Z}_{[S]}$ with a subset of each
  stalk.
\end{lemma}
A section of $\mathcal{Z}$ over a base $X$ is said to be \emph{nodal}
if it restricts to a node in some stalk.  Lemma
\ref{central-stalk-lemma} implies that the nodal sections over
$\fB(S)$ are precisely those which restrict to nodes at the origin
$[S\to S]$.

\begin{lemma}\label{marking-and-collapse-lemma}
Given a point $[F \stackrel{\alpha}{\to} S] \in \fB(S)$, the marking
$\alpha$ collapses to nodes precisely those curves in $F$ which are
the restrictions of nodal sections over $\fB(S)$
\end{lemma}
\begin{proof}
  By Lemma \ref{central-stalk-lemma}, a curve in $F$ is the preimage
  of a node in $S$ if and only if it is the restriction of a nodal
  section over $\fB(S)$.
\end{proof}

Let $T$ be a stable nodal surface and let $\alpha_*\co T\to S$ be a
degeneration which collapses a single curve in $T$ to a node $n \in
S$.  The node $n$ determines a nodal section over $\fB(S)$, and the
maximal subset over which this section is \emph{not} nodal is
precisely the image of the change-of-marking $\alpha_*\co \fB(T) \to
\fB(S)$.

\begin{lemma}\label{sections-identification}
  Suppose $X$ is simply connected and $\widetilde{\sigma}\co X \to
  \fB(S)$ is a maximal lift.  Then the map $\widetilde{\sigma}^*\co
  \mathcal{Z}(\fB(S)) \hookrightarrow \mathcal{Z}(X)$ restricts to a
  bijection between nodal sections.
\end{lemma}
\begin{proof}
  Clearly every nodal section over $X$ is the pullback of a nodal
  section over $\fB(S)$.  Conversely, if there exists a nodal section
  over $\fB(S)$ which pulls backs to a non-nodal section over $X$ then
  $X$ lies in the image of a change-of-marking $\alpha_*\co \fB(T) \to
  \fB(S)$ for some $T$ with strictly fewer nodes than $S$.  Since $X$
  is simply connected and the change-of-marking maps are local
  homeomorphisms, $X$ lifts further, which contradicts the maximality
  hypothesis.
\end{proof}

\begin{proof}[Proof of Lemma \ref{initial-lifting-lemma}]
  Choose a point $x\in X$ and let $F_x$ denote the fiber over $x$.
  Consider the commutative diagram
  \begin{equation*}
    \begin{diagram}
      \node[2]{\mathcal{Z}(X)} \arrow{s,J} \\
      \node{\mathcal{Z}(\fB(S))} \arrow{ne,t,J}{\sigma_1^*}
      \arrow{e,J} \node{\mathcal{Z}_{[F_x]}}
    \end{diagram}
  \end{equation*}
  where the vertical and horizontal arrows are induced by restriction
  to the stalk at $x$ (which is identified with the stalk at
  $\sigma_1(x)$.  By Lemma \ref{sections-identification}, a curve in
  $F_x$ is the restriction of a nodal section over $X$ if and only if
  it is the restriction of a nodal section over $\fB(S)$.  By Lemma
  \ref{marking-and-collapse-lemma}, $S$ is topologically obtained from
  $F_x$ by collapsing those curves which are the restrictions of nodal
  sections over $\fB(S)$ (equivalently, nodal sections over $X$).  By
  the same reasoning, $T$ is topologically obtained from $F_x$ by
  collapsing the same set of curves.  Hence $S$ and $T$ are abstractly
  homeomorphic.  Finally, since the image of $\fB(S)$ in $\bMan$ is
  isomorphic to the quotient stack $[\fB(S) / \MCG(S)]$, it follows
  that any two lifts to $\fB(S)$ are related by a unique change of
  marking.
\end{proof}

\section{Proof of Lemma \ref{inclusion-htpy-equiv-lemma}}\label{proof-of-main-lemma}

We shall now prove Lemma \ref{inclusion-htpy-equiv-lemma}.  It will
follow from Lemma \ref{initial-lifting-lemma} together with
Waldhausen's Theorem $A'$ \cite[p. 165]{Waldhausen-manifolds}, which
is a simplicial version of Quillen's Theorem A.  We first recall
Waldhausen's Theorem $A'$. Suppose $F\co A_\bullet \to B_\bullet$ is a
functor of simplicial categories.  Given an object $\sigma \in \Obj
B_n$ of simplicial degree $n$, the simplicial fiber category
$(F/\sigma)_\bullet$ is given in degree $k$ by
\[
(F/\sigma)_k := \coprod_{u\co \underline{k} \to \underline{n}} F_k/u^* \sigma,
\]
where the disjoint union is taken over all monotone maps from $\{0,
\ldots k\}$ to $\{0 \ldots n\}$.  The theorem states that if each of
these simplicial fiber categories has contractible classifying space
then $F$ is a homotopy equivalence of classifying spaces.

\begin{proof}[Proof of Lemma \ref{inclusion-htpy-equiv-lemma}]
By taking the total singular simplicial set one has an inclusion of
simplicial categories
\[
j\co S_\bullet\left( \CL_{g,n} \int \widehat{\fB} \right) \rightarrow
S_\bullet \fB.
\]
We will apply Waldhausen's Theorem $A'$ to the simplicial functor $j$.
Fix an object $\phi\co \Delta^n \to \Obj \: \fB$ of $S_n \fB$.  The
image of $\phi$ lands in $\fB(S)$ for some stable nodal surface $S$.
In simplicial degree $k$ the simplicial fiber category
$(j/\phi)_\bullet$ is a disjoint union of ordinary fiber categories of
the form $j_k/\sigma$ for various objects $\sigma$ of simplicial
degree $k$.  Suppose that each of these categories is contractible.
Then collapsing them to points maps the simplicial fiber category
$(j/\phi)_\bullet$ by a levelwise homotopy equivalence to the standard
simplicial model for the $n$-simplex given in degree $k$ by
$\coprod_{u\co \underline{k} \to \underline{n}} *$.  The geometric
realization of this map is thus a homotopy equivalence
$|B(j/\phi)_\bullet| \to
\Delta^n \simeq *,$ and so Waldhausen's Theorem A' yields the result.

It thus suffices to show that each category $j_n/\sigma$ has an
initial object, where $\sigma\co \Delta^n \to \fB(S) \subset \Obj \: \fB$.
Let $\overline{\sigma}\co \Delta^n \to \bMan$ denote the composition of
$\sigma$ with the projection of $\fB(S)$ down to $\bMan$.  Explicitly,
an object of the category $j_n/\sigma$ is a lift (up to a
specified 2-morphism) of $\overline{\sigma}$ to some chart $\fB(T)$;
i.e. a 2-commutative diagram
\[
\xymatrix{
 & \fB(T) \ar[d] \\
\Delta^n \ar[r]_{\overline{\sigma}} \ar[ur]^{\tau}_{\Downarrow \theta} & \bMan.
}
\]
A morphism $(\tau_1,\theta_1) \to (\tau_2,\theta_2)$ is an isotopy
class of degenerations $\alpha\co T_1 \to T_2$ such that the induced
2-morphism $\Phi(\alpha)$
\[
\begin{diagram}
\node[2]{\fB(T_1)} \arrow{se} \\
\node{\Delta^n} \arrow{ne} \arrow{se} \node{\Downarrow \Phi(\alpha)} \node{\bMan} \\
\node[2]{\fB(T_2)} \arrow{ne}
\end{diagram}
\]
satisfies $\theta_2 \circ \Phi(\alpha) = \theta_1$.  This is
equivalent to saying that $\alpha_* \circ \tau_1 = \tau_2$.  Since
$j_n$ is the inclusion of a subcategory into a groupoid, there is at
most one arrow between any two objects of $j_n / \sigma$.  We are thus
reduced to showing that there is an object $(\sigma_0, \theta_0)$
which maps to all other objects of $j_n / \sigma$.  Every lift of
$\overline{\sigma}$ to some $\fB(T)$ lifts further to a maximal lift,
and Lemma \ref{initial-lifting-lemma} says that a maximal lift of
$\overline{\sigma}$ is unique up to isomorphism.  A maximal lift
therefore provides the desired initial object.
\end{proof}

\section{Stratifications}\label{Stratifications}

The strata of $\bMan$ are indexed by \emph{stable graphs} with $n$
external legs; equivalently the strata are indexed by diffeomorphism
types of stable nodal surfaces of genus $g$ with $n$ labeled points.
A stable nodal surface $T$ corresponds to an open stratum $R_T\bMan$
which is the locus of all stable nodal Riemann surfaces $F$ which are
topologically diffeomorphic to $T$.  The closure $\overline{R}_T\bMan$ is
the locus of all Riemann surfaces $F$ for which $T$ admits a
degeneration onto $F$.  This stratification gives a corresponding
stratification of the spaces $\fB(S)$, and so there are atlases:
\begin{align*}
\coprod_S R_T \fB(S) &\to R_T \bMan, \\
\coprod_S \overline{R}_T \fB(S) &\to \overline{R}_T \bMan
\end{align*}
which give rise to subgroupoids of the Bers groupoid $\fB$.  The
Fenchel-Nielsen coordinates show that $\overline{R}_T \fB(S)$ is
homeomorphic to a proper ball in $\fB(S)$ of codimension equal to the
number of nodes of $S$ minus the number of nodes of $T$.  In
particular, $\overline{R}_T \fB(S)$ is contractible.

The stratification of $\CL_{g,n}$ is as follows:
\begin{align*}
R_T \CL_{g,n} & = \mbox{full subcategory on the object $T$} = \MCG(T)\\
\overline{R}_T \CL_{g,n} & = \{\mbox{full subcategory on $S$
such that $T$ admits a degeneration onto $S$}\}
\end{align*}

The proof of Theorem \ref{main-theorem} (along with the proofs of all
propositions and lemmas it employs) remains valid upon inserting
$\overline{R}_T$ in front of all occurences of the symbols
$\CL_{g,n}$, $\bMan$, and $\fB$.  Thus the homotopy equivalence $\bMan \simeq
B\CL_{g,n}$ restricts to an equivalence of each closed stratum.

To see that it restricts to a homotopy equivalence on each open
stratum, one uses the fact that each open stratum $R_T \bMan$ is the
stack quotient of a finite group acting on a product of uncompactified
moduli spaces, and $R_T \CL_{g,n}$ is the homotopy quotient of the
same finite group acting on the corresponding product of classifying
spaces of mapping class groups.  Thus $R_T \bMan \simeq B(R_T
\CL_{g,n})$ follows from the equivalence $\M^{an} \simeq
B\MCG$ discussed in the introduction, since the homotopy type
of the stack quotient is the homotopy quotient.

\bibliographystyle{amsalpha}
\bibliography{biblio}
\end{document}